\documentclass[onefignum,onetabnum]{siamonline171218}

\usepackage{lipsum}
\usepackage{amsfonts}
\usepackage{amssymb}
\usepackage{graphicx}
\usepackage{epstopdf}
\usepackage{algorithmic}
\ifpdf
  \DeclareGraphicsExtensions{.eps,.pdf,.png,.jpg}
\else
  \DeclareGraphicsExtensions{.eps}
\fi

\usepackage{enumitem}
\setlist[enumerate]{leftmargin=.5in}
\setlist[itemize]{leftmargin=.5in}

\newsiamremark{remark}{Remark}
\newsiamremark{hypothesis}{Hypothesis}
\crefname{hypothesis}{Hypothesis}{Hypotheses}
\newsiamthm{claim}{Claim}

\headers{CDMD}{Azencot et al.}

\title{Consistent Dynamic Mode Decomposition\thanks{Submitted to the editors DATE.
\funding{This work was supported by the European Union\textsc{\char13}s Horizon 2020 research and innovation programme under the Marie Sk{\l}odowska-Curie grant agreement No. 793800, a Zuckerman STEM Leadership Postdoctoral Fellowship, NSF Grant DMS-1720237, ONR Grant N000141712162, NSF grant DMS-1737770 and the City of Los Angeles, Gang Reduction Youth Development (GRYD) Analysis Program.}}}

\author{Omri Azencot \and Wotao Yin \and Andrea Bertozzi
\thanks{Department of Mathematics, University of California, Los Angeles CA 90095, (\email{azencot@math.ucla.edu}, \email{wotaoyin@math.ucla.edu}, \email{bertozzi@math.ucla.edu}).}}

\usepackage{amsopn}
\DeclareMathOperator*{\argmin}{arg\,min}
\DeclareMathOperator*{\minz}{minimize}
\DeclareMathOperator{\Tr}{Tr}

\DeclareMathOperator{\im}{Im}

\DeclareMathOperator{\vecm}{{\mathrm{vec}}}

\newcommand{\mca}[1]{\mathcal{#1}}
\newcommand{\mbb}[1]{\mathbb{#1}}

\newcommand{\code}[1]{\texttt{#1}}

\ifpdf
\hypersetup{
  pdftitle={},
  pdfauthor={}
}
\fi

\begin{document}

\maketitle

\begin{abstract}
We propose a new method for computing Dynamic Mode Decomposition (DMD) evolution matrices, which we use to analyze dynamical systems. Unlike the majority of existing methods, our approach is based on a variational formulation consisting of data alignment penalty terms and constitutive orthogonality constraints. Our method does not make any assumptions on the structure of the data or their size, and thus it is applicable to a wide range of problems including non-linear scenarios or extremely small observation sets. In addition, our technique is robust to noise that is independent of the dynamics and it does not require input data to be sequential. Our key idea is to introduce a regularization term for the forward and \emph{backward} dynamics. The obtained minimization problem is solved efficiently using the Alternating Method of Multipliers (ADMM) which requires two Sylvester equation solves per iteration. Our numerical scheme converges empirically and is similar to a \emph{provably} convergent ADMM scheme. We compare our approach to various state-of-the-art methods on several benchmark dynamical systems.

\end{abstract}



\section{Introduction}

Over the last few years, data-driven approaches became prevalent in analyzing dynamical systems~\cite{kutz2016dynamic}. In the common scenario, a collection of system \emph{observations} is provided and a \emph{linear} object that encodes the dynamics is generated based solely on the data. These data-driven approaches are advantageous in that they make minimal assumptions on the governing equations of the system, and in particular, these techniques are applicable even to non-linear dynamics. In this context, Dynamic Mode Decomposition (DMD)~\cite{schmid2010dynamic} methods gained a lot of attention lately, in part due to their computational efficiency as well as their analysis capabilities of the system at hand. DMD-based methods were successfully applied to various flows including detonation waves, cavity flows and jets~\cite{massa2012dynamic,seena2011dynamic,schmid2011applications}. In short, DMD computes a matrix whose spectrum, represented by the eigenvalues and eigenvectors, provides meaningful information such as growth and decay rates of the system or dominant coherent structures in the flow. The goal of this paper is to propose a new method for computing DMD matrices that is based on interpreting the problem in a variational form, taking into account the forward and backward dynamics and solving it efficiently via splitting.

Developing data-driven methodologies for the analysis of non-linear dynamical systems is an active research domain with DMD being one of its main avenues. In particular, DMD was recently generalized and extended in several works having the objective of alleviating some of the shortcomings in the original technique. For instance, a limiting assumption in~\cite{rowley2009spectral,schmid2010dynamic} requires that the data is given in a sequential form, namely, the input snapshots represent an equally spaced time series of observations. In Tu et al.~\cite{tu2013dynamic} and other works, this limitation is relaxed and pairs of equispaced observations are used instead, whereas in~\cite{gueniat2015dynamic,leroux2016dynamic,askham2018variable}, no assumption is made on the regularity of the temporal sampling. Another drawback of several DMD methods is the bias they exhibit in the presence of noise and whether the noise interacts with the dynamics~\cite{bagheri2014effects} or not~\cite{dawson2016characterizing}. To address this challenge, variants of DMD were proposed in the literature based on solving jointly for the basis and the evolution operator~\cite{wynn2013optimal}, formulating the problem as a total least squares minimization~\cite{hemati2017biasing}, and fitting an exponential model~\cite{askham2018variable}. Other methods cope with noise by utilizing Kalman filters~\cite{nonomura2018dynamic,nonomura2019extended}, adapting DMD to online data~\cite{hemati2014dynamic,hemati2016improving}, and developing a Rayleigh--Ritz modal decomposition~\cite{drmac2018data}, among other approaches~\cite{dawson2016characterizing}. Under this classification, our method is applicable to non-sequential data and it performs extremely well when sensor noise 
corrupts the data, as we show in Section~\ref{sec:results}. 

Perhaps closest to our approach is the work of Dawson et al.~\cite{dawson2016characterizing} where the idea of making DMD more robust to noise by considering the forward and \emph{backward} evolution is investigated. More specifically, in forward-backward DMD (fbDMD)~\cite{dawson2016characterizing}, the DMD matrix is estimated via the square root of the product of the forward model with the inverse of the backward DMD matrix. The backward estimate is generated by switching the ``before'' and ``after'' roles of the snapshots. Our machinery is based on the same observation of exploiting the forward and backward dynamics, but in a completely different way. Inspired by ideas from \emph{Computer Graphics}~\cite{ovsjanikov2012functional,huang2014functional}, we formulate the task of computing the DMD matrix in a variational form that includes penalties for both directions. The obtained minimization is unfortunately highly non-linear and non-convex, and thus we introduce an auxiliary variable that represents the backward dynamics, arriving at an optimization problem with quadratic objective terms and bilinear constraints. This problem can be solved efficiently using splitting techniques such as the Alternating Direction Method of Multipliers (ADMM)~\cite{boyd2011distributed}. The obtained scheme is iterative, where at each step we solve two Sylvester equations and perform a trivial update. In addition, we show that our problem can be modified such that a provably convergent scheme can be devised. Overall, we obtain an efficient algorithm that exhibits fast convergence rates in practice and provides improved estimates of various properties of the dynamical system.
 
The rest of the paper is organized as follows. In~\Cref{sec:background} we provide background details related to dynamic mode decomposition techniques and the alternating method of multipliers. \Cref{sec:cdmd} details our approach for generating consistent DMD evolution matrices where we derive the variational formulation, and we propose an effective ADMM splitting scheme to solve it in practice. In~\Cref{sec:converging_cdmd}, we prove that the problem we consider can be changed so that it admits an ADMM-type algorithm which is provably converging. \Cref{sec:results} provides a quantitative and qualitative evaluation of our method with respect to several DMD algorithms. \Cref{sec:conclusions} concludes our work, discusses limitations, and offers a few potential directions for future work.

\section{Background}
\label{sec:background}

In what follows, we briefly present the most relevant details regarding DMD algorithms. We refer to~\cite{kutz2016dynamic} for a more comprehensive text on the recent developments and applications of DMD-based techniques. In addition, we describe the essential components of ADMM and their link to our work, where we point to the paper by Boyd et al.~\cite{boyd2011distributed} for additional information.

\begin{figure*}[t]
 \centering
 \includegraphics[width=\linewidth]{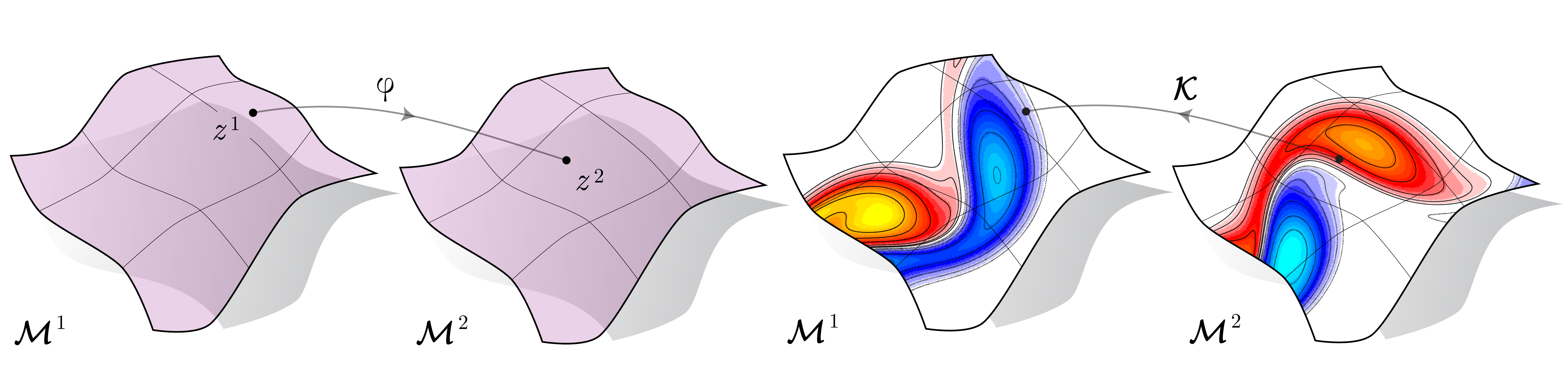}
 \vspace{-4mm}
 \caption{The dynamical system $\varphi$ maps from the manifold $\mca{M}$ at time $1$ to time $2$ (left), whereas the associated Koopman operator maps between scalar functions defined on $\mca{M}$ (right).}
 \label{fig:koopman}
\end{figure*}

\subsection{DMD} 
\label{subsec:dmd}

Dynamic Mode Decomposition (DMD) emerged in the fluid dynamics field~\cite{schmid2010dynamic} as a data driven approach for analyzing a dynamical system based on observational data. DMD is strongly related to Koopman theory~\cite{koopman1931hamiltonian}, where a non-linear dynamical system $\varphi$ acting on a finite-dimensional manifold $\mca{M}$ is encoded using an infinite-dimensional linear Koopman operator $\mca{K}$. In this context, DMD can be viewed as a practical approach to produce a matrix $A$ whose spectrum approximates the spectrum of the operator $\mca{K}$. Thus, $A$ is an informative object and its dominant eigenvalues and eigenvectors are directly linked to dynamical features of the system such as growth, decay, frequency and flow modes. These results encourage the community to investigate DMD as an effective tool for analyzing various linear and nonlinear dynamical systems~\cite{kutz2016dynamic}.

A common scenario, considered in several DMD-based techniques, is to assume to be given a set of temporally related pairs of observations $\tilde{x}_j$ and $\tilde{y}_j, j = 1,2,..,n$, such that
\begin{align}
\tilde{y}_j( z ) = \tilde{x}_j( \varphi (z) ) \ ,  
\end{align}
where $z \in \mca{M}$, the dynamical system is $\varphi : \mca{M} \rightarrow \mca{M}$, and $\tilde{x}_j,\, \tilde{y}_j : \mca{M} \rightarrow \mbb{R}$. Namely, if $\tilde{x}_j$ represents some quantity at time $t$, then $\tilde{y}_j$ measures the same quantity at a later time $t+\Delta t$, as it changes due to the dynamics $\varphi$, see Fig.~\ref{fig:koopman} for an illustration of this setup. Examples of the input observations could be the spatial coordinates~\cite{dawson2016characterizing} or the scalar vorticity~\cite{tu2013dynamic}, among other system-related data. The time series of observations $\{ \tilde{x}_j \}_{j=1}^n$ and $\{ \tilde{y}_j \}_{j=1}^n$ is used to construct matrices $\tilde{X}$ and $\tilde{Y}$ such that
\begin{align}
\tilde{X} = [\tilde{x}_1 \; \tilde{x}_2 \; ... \; \tilde{x}_n] \in \mbb{R}^{m \times n}, \quad \tilde{Y} = [\tilde{y}_1 \; \tilde{y}_2 \; ... \; \tilde{y}_n] \in \mbb{R}^{m \times n} \ ,    
\end{align}
where the manifold $\mca{M}$ is of dimension $|\mca{M}|=m$. We note that our data is equispaced in time, i.e., $\Delta t$ is the same for every $j$, as is commonly assumed in the DMD literature, although other scenarios were considered, e.g.,~\cite{tu2014spectral}. Using the above notation, the goal of many DMD algorithms is to find a matrix $\tilde{A} \in \mbb{R}^{m \times m}$ such that $\tilde{A} \tilde{X} = \tilde{Y}$.

In practice, solving directly for $\tilde{A}$ could be challenging, especially when $m$ is extremely large or when $m>n$, leading to an underdetermined system. One way to mitigate these difficulties is to reduce the spatial dimension of the input data. Many dimensionality reduction techniques have been developed in recent years, where the Proper Orthogonal Decomposition (POD)~\cite{berkooz1993proper} is typically chosen mostly due to its algorithmic simplicity and computational efficiency. One of the outputs of POD is a set of $r$ orthogonal modes $B \in \mbb{C}^{m \times r}$ such that the linear subspace spanned by $B$ approximates $\mbb{R}^m$ well enough. From now on, we denote by $X$ and $Y$ the projection of $\tilde{X}$ and $\tilde{Y}$ onto the first $r$ POD modes. Formally,
\begin{align}
X = B^* \tilde{X} \in \mbb{R}^{r \times n} \ , \quad Y = B^* \tilde{Y} \in \mbb{R}^{r \times n} \ ,
\end{align}
where $B^*$ is the conjugate transpose of $B$. To compute the matrix $B$, we facilitate the Singular Value Decomposition (SVD) to obtain the expression $\tilde{X} = \tilde{U} \tilde{S} \tilde{V}^*$ and $B = \tilde{U}_r$, i.e., the first $r$ left singular vectors that correspond to the dominant $r$ singular values. In this reduced form, the problem of DMD is to solve the equation $AX=Y$, for which the least squares solution is analytically given by $A=YX^{+}$, where $X^{+}$ is the Moore--Penrose pseudoinverse of $X$. The DMD algorithms we will present next can be thought of as various approaches for approximating such $A$ matrices.

The following \cref{alg:exact_dmd} was introduced in Tu et al.~\cite{tu2013dynamic} and is known as the \code{Exact DMD} method. While this approach is not one of the original DMD techniques as was proposed in \cite{rowley2009spectral,schmid2010dynamic}, it is a close variant of these methods and it serves as the baseline algorithm for many extensions and comparisons in the DMD literature. We note that in Step~($3$), instead of taking the pseudoinverse $X^{+}$, the authors took its projection onto the first $r$ modes. Indeed, we have that
\[
X^{+} = \left( B^* \tilde{X} \right)^{+} = \left( \tilde{U}_r^* \tilde{U} \tilde{S} \tilde{V}^* \right)^{+} \approx \left(\tilde{U}_r^* \tilde{U_r} \tilde{S_r} \tilde{V_r}^* \right)^{+} = \tilde{V}_r \tilde{S}_r^{-1} \ .
\]
Also, Step~($4$) involves the eigendecomposition (EIG) of $A$, typically yielding a complex-valued spectrum since $A$ does not exhibit a special structure in general. Finally, in many DMD-based algorithms, Steps~($1-2$) and ($4-5$) are shared, whereas Step~($3$) is different. This is also the case in our~\cref{alg:cdmd} where the main change is the way we construct the matrix $A$.

\begin{algorithm}
\caption{Exact Dynamic Mode Decomposition (\code{Exact DMD})}
\label{alg:exact_dmd}
\begin{algorithmic}[1]
\STATE{Input matrices $\tilde{X},\tilde{Y} \in \mbb{R}^{m\times n}$ and a scalar $r\in \mbb{R}$}
\vspace{2mm}
\STATE{Compute the SVD of $\tilde{X}=\tilde{U} \tilde{S} \tilde{V}^*$, and generate $X = \tilde{U}_r^* \tilde{X}, Y = \tilde{U}_r^* \tilde{Y}$}
\vspace{2mm}
\STATE{Denote $A = Y \tilde{V}_r \tilde{S}_r^{-1}$}
\vspace{2mm}
\STATE{Compute the EIG of $A$, with $A v_j = \lambda_j v_j$, where $v_j \in \mbb{C}^r, \; \lambda_j \in \mbb{C}$}
\vspace{2mm}
\STATE{The DMD spectrum is defined as the set of eigenvalues $\lambda_j$, and vectors $\psi_j = \lambda_j^{-1} \tilde{Y} \tilde{V}_r \tilde{S}_r^{-1} v_j$}

\end{algorithmic}
\end{algorithm}

\subsection{Regularizing DMD}
\label{subsec:reg_dmd}

In many scenarios, the time sequence of data is generated using sensory devices. For example, Schmid et al.~\cite{schmid2011applications} applied DMD to snapshots of a helium jet, collected using particle-image-velocimetry (PIV) measurements. Naturally, in these settings, the observations are assumed to be corrupted with various types of noise. The existence of process or sensor noise results in a certain bias in traditional DMD algorithms such as \code{Exact DMD}, as was recently shown in~\cite{dawson2016characterizing,hemati2017biasing,askham2018variable}. To address these shortcomings, several extensions to DMD were recently proposed in the literature. From an optimization standpoint, these modified DMD methods as well as our approach can be viewed as \emph{regularizing} the original minimization problem, introducing algorithms that are more robust in the presence of noise. In our discussion here, we focus on the methods \code{fbDMD}~\cite{dawson2016characterizing}, \code{tlsDMD}~\cite{hemati2017biasing} and \code{Optimized DMD}~\cite{askham2018variable}.

The main idea behind the forward-backward DMD (\code{fbDMD}) technique is to take into account the forward dynamics, i.e., transforming $X$ into $Y$, as well as the backward system where $Y$ is mapped to $X$. The motivation is that by considering both directions, much of the bias to noise can be eliminated. In fact, we build on the exact same observation, however, we arrive at a completely different method. The algorithm \code{fbDMD} follows the same steps of~\cref{alg:exact_dmd}, except for the matrix construction which is given by
\begin{align}
A = \left( A_f A_b^{-1} \right)^{1/2} \ ,
\end{align}
where $A_f = \tilde{U}_X^* \tilde{Y} \tilde{V}_X \tilde{S}_X^{-1}$ is the forward estimate, and $A_b = \tilde{U}_Y^* \tilde{X} \tilde{V}_Y \tilde{S}_Y^{-1}$ is the backward one. Notice that the SVD of both $\tilde{X} = \tilde{U}_X \tilde{S}_X \tilde{V}_X^*$ and $\tilde{Y} = \tilde{U}_Y \tilde{S}_Y \tilde{V}_Y^*$ are used. Assuming that efficient routines for computing the square root of a matrix such as \code{sqrtm} of MATLAB are available, the time complexity for this algorithm is $\mca{O} \left( \min\{mn^2,m^2n\}+r^3 \right)$, and thus it is governed by the SVD part as we typically have $r \ll m,n$.

In a different paper~\cite{hemati2017biasing}, the authors propose another algorithm known as the total least squares DMD (\code{tlsDMD}). Intuitively, this approach tries to symmetrize the way noise is being handled so that it assumes noise polluted both $X$ and $Y$, whereas other methods implicitly account only for noise in $Y$. Similarly to the latter algorithm, \code{tlsDMD} provides an alternative definition for the $A$ matrix. Specifically,
\begin{align} \label{eq:tlsDMD}
A = U_{br} \, U_{tr}^{-1} \ , \quad \text{with} \; \begin{pmatrix} X \\ Y \end{pmatrix} = U S V^* \; \text{and} \; U = \begin{pmatrix} U_{tr} & \mca{U}_{tr} \\ U_{br} & \mca{U}_{br} \end{pmatrix} \ .
\end{align}
Namely, the projected observations $X$ and $Y$ are combined into a matrix of size $2r \times n$, whose $r$ dominant left singular vectors are used to compute $A$. The matrix $U_{tr} \in \mbb{C}^{r \times r}$ encodes the top left part of $U$ and $U_{br} \in \mbb{C}^{r \times r}$ represents the bottom left part of $U$. The scalar $r$ satisfies $r < n/2$ in this method. Overall, the computational requirements of \code{tlsDMD} are on the order of $\mca{O} \left( \min\{mn^2,m^2n\}+r^3 \right)$.

Finally, a recent development for computing DMD matrices was introduced in~\cite{askham2018variable} resulting in the \code{Optimized DMD} method. Essentially, the authors formulate DMD as a non-linear least squares minimization problem. To this end, the ensemble of observations is put together, e.g., $Z = \begin{pmatrix} X & Y \end{pmatrix} \in \mbb{R}^{m \times 2n}$, and the goal is to fit $Z$ with a linear combination of non-linear functions $\Phi \in \mbb{R}^{2n \times l}$. In practice, $\Phi$ is taken from a family of exponential functions such as $\Phi(\alpha,t)_j = \exp(\alpha_j t)$, where the set of parameters $\alpha \in \mbb{C}^k$ is unknown. The optimization problem takes the form of
\begin{align} \label{eq:optimizedDMD}
\minz_{\alpha,B} \quad |Z^T - \Phi(\alpha) B|_F^2 \ ,
\end{align}
where $B \in \mbb{C}^{l \times m}$ is the set of unknown coefficients which determine the linear superposition of non-linear functions from $\Phi$. Observing that $B$ can be eliminated from the optimization, problem~\cref{eq:optimizedDMD} may be efficiently solved using the \emph{variable projection method}~\cite{golub1973differentiation}. We note that the DMD spectrum and the matrix $A$ could be constructed using the computed outputs $\Phi$ and $B$, and we refer to~\cite{askham2018variable} for further details.

\subsection{ADMM}

The \emph{Alternating Direction Method of Multipliers} (ADMM) is a numerical optimization approach for efficiently solving separable objective functions. ADMM was first introduced in 1970's in \cite{GlowinskiMarroco1975_LapproximationPar,GabayMercier1976_dual}, recently popularized by \cite{GoldsteinOsher2009_split, boyd2011distributed}, and generalized for nonconvex optimization in \cite{WangYinZeng2015_global,gao2018admm}). A general scenario for which ADMM is effective involves the following minimization problem,
\begin{align} \label{eq:lin_sep_pbm}
\minz_{x, z} \quad f(x) + g(z) \ , \quad \text{s.t.} \quad Ax + Bz = c \ ,
\end{align}
where $f(x):\mbb{R}^n \rightarrow \mbb{R}$ and $g(z):\mbb{R}^m \rightarrow \mbb{R}$ are convex functions, the linear constraints include matrices $A \in \mbb{R}^{p \times n}$, $B\in\mbb{R}^{p \times m}$ and a vector $c \in \mbb{R}^p$. To solve~\cref{eq:lin_sep_pbm}, we define the following augmented Lagrangian,
\begin{align}
\mca{L}_\rho(x,z,y) = f(x) + g(z) + y^T(Ax+Bz-c) + \frac{\rho}{2} |Ax+Bz-c|_2^2 \ .
\end{align}
ADMM exploits the fact that $\mca{L}_\rho$ can be decomposed with respect to the variables $x$ and $z$, leading to a numerical splitting scheme consisting of the iterations
\begin{align} \label{eq:lin_admm} \begin{split}
x^{k+1} &= \argmin \mca{L}_\rho(x,z^k,y^k) \\
z^{k+1} &= \argmin \mca{L}_\rho(x^{k+1},z,y^k) \\
y^{k+1} &= y^k + \rho ( Ax^{k+1}+B z^{k+1}-c ) \ ,
\end{split} \end{align}
where $\rho>0$ is the penalty parameter in the augmented Lagrangian. The advantage of utilizing ADMM is twofold, solving alternately for $x$ and $z$ typically involves simpler minimization problems compared to a joint optimization, and convergence results require mild assumptions.

It is often useful to facilitate a change of variables and to define a \emph{scaled} version for the dual variable $y$, denoted by $\rho u = y$. This choice significantly reduces the length of formulas, and thus we will opt for this version throughout the paper. We denote by $r(x,z) = Ax+Bz-c$, and we re-write the scaled augmented Lagrangian in terms of $u$,
\[
\mca{L}_\rho(x,z,u) = f(x) + g(z) + \frac{\rho}{2} |r + u|_2^2 - \frac{\rho}{2} |u|_2^2 \ .
\]
The associated splitting scheme is similar in the $x$ and $z$ updates where we replace $y^k$ with $u^k$ in~\cref{eq:lin_admm}, whereas for the $u$ update we have $u^{k+1} = u^k + r(x^{k+1},z^{k+1})$.

\section{Consistent Dynamic Mode Decomposition}
\label{sec:cdmd}

In this section we describe our main algorithm for computing an approximation of the DMD operator that is associated with some known dynamical observations. The key observation in our approach is the consideration of the forward \emph{and} backward dynamics within the same framework. In this context, we propose a variational formulation of the problem where we simultaneously solve for the forward and backward DMD operators. Unfortunately, the formulation we arrive at is highly non-linear and non-convex, and thus challenging to solve in practice. Our main contribution is an effective splitting numerical scheme which is efficient yet easy to code.

\subsection{Forward and backward dynamics} 

Let the two matrices $X, Y \in \mbb{R}^{r\times n}$ represent our POD-projected data such that each column in $X$ is associated with the corresponding column in $Y$ under the dynamics (see~\Cref{subsec:dmd}). Several Dynamic Mode Decomposition (DMD) algorithms study the forward dynamics, i.e., find $A$ such that $A X \approx Y$. We advocate the consideration of the backward dynamics, namely, we also want that $A^{-1} Y \approx X$. This idea was previously explored in~\cite[Section 2.4]{dawson2016characterizing}, where the authors proposed the \code{fbDMD} algorithm which takes into account both directions. However, there are a few key differences between our approach and theirs, as we detail below. Formally, we consider the following variational problem,
\begin{align} \label{eq:fb_pbm}
    \minz_{A} \quad \frac{1}{2} \left|AX-Y\right|_F^2 + \frac{1}{2} \left|X-A^{-1}Y\right|_F^2 \ ,
\end{align}
where $|\cdot|_F$ is the Frobenius norm. We note that if $A$ is orthogonal, i.e., $A^{-1} = A^T$, then the above addends are equal, however in the general case we have
\begin{align*}
    \left|AX-Y\right|_F^2 &= \Tr(X^TA^TAX - 2 X^TA^TY + Y^TY) \\
    &\neq \Tr(X^TX - 2 X^TA^{-1}Y + Y^TA^{-T}A^{-1}Y) = \left|X-A^{-1}Y\right|_F^2 \ .
\end{align*}

\begin{figure*}[t]
 \centering
 \includegraphics[width=\linewidth]{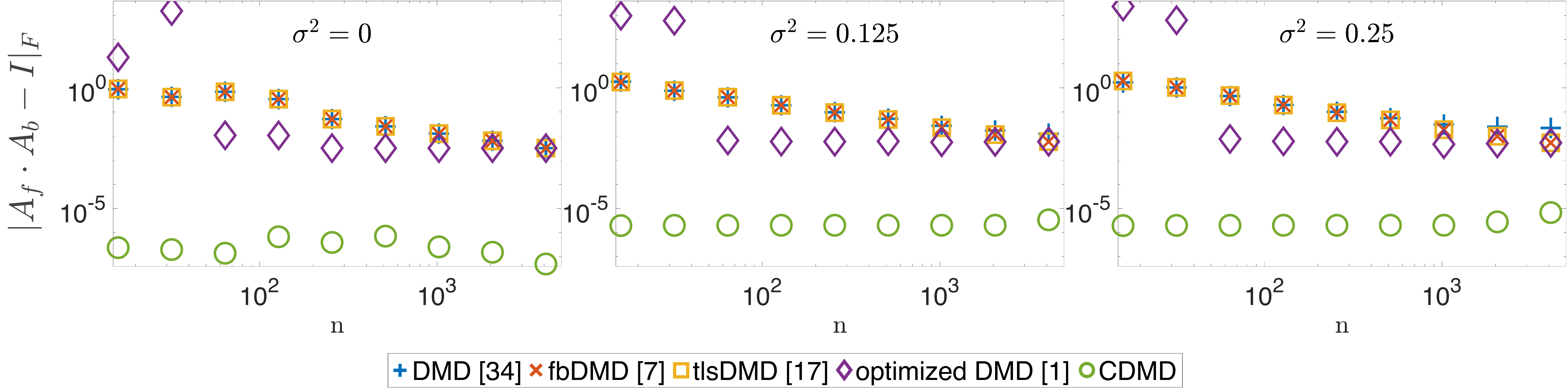}
 \vspace{-8mm}
 \caption{Our method penalizes the obtained inconsistency of composing the forward and backward evolution operators, $A, B$. In practice, employing consistency constraints in our optimization regularizes the problem significantly, yielding robust estimates in the presence of noise as we show in \Cref{sec:results}. In the plots above, we demonstrate the consistency error for the sine example~\cref{eq:nonlin_sine} as achieved by various DMD techniques. Our method yields extremely low error rates, whereas \code{Exact DMD}, \code{fbDMD} and \code{tlsDMD} obtain higher rates that decrease as the number of observations $n$ increases. Finally, \code{optimized DMD} generates the second to best consistency estimates with the exception of low number of observations, where their error rates are the highest. }
 \label{fig:sine_consistency}
\end{figure*}

\subsection{Change of variables} 

The optimization problem~\cref{eq:fb_pbm} is highly non-linear and non-convex due to the $A^{-1}$ term. Therefore, instead of directly solving this challenging problem, we introduce the auxiliary variable $B=A^{-1}$, and we re-formulate to arrive at, 
\begin{align} \label{eq:fb_pbm2}
    \minz_{A,B} \quad \frac{1}{2} \left|AX-Y\right|_F^2 + \frac{1}{2} \left|X-BY\right|_F^2 \ , \quad \text{s.t.} \quad AB = I, BA = I \ ,
\end{align}
where the constitutive constraints $AB=I$ and $BA=I$ guarantee that minimizers of~\cref{eq:fb_pbm2} are inverse of each other. From an optimization point of a view, if one of the constraints is satisfied then the second constraint holds as well. However, in practice, adding both constraints is a reasonable choice as they symmetrize the approximate invertible relations of $A$ and $B$. We refer to the above re-formulation as the \emph{Consistent Dynamic Mode Decomposition (CDMD)} problem. To motivate our methodology, we quantify the consistency error $|AB-I|_F$ obtained by several existing methods including ours, and we plot the results in \Cref{fig:sine_consistency}. Indeed, our technique is highly consistent compared to the other approaches, almost independently of the number of observations $n$. We note that when the consistency error is large, it may hint of overfitting to data, since the forward and backward estimations represent systems that are far from being inverse of each other. 

The CDMD functional~\cref{eq:fb_pbm2} appeared previously in Computer Graphics applications where a discrete map between two dimensional surfaces is being sought. Namely, given two geometric shapes such as two different poses of the same person, the goal is to determine where each point on one shape is mapped to its corresponding point on the second shape. DMD operators (also known as \emph{functional maps}~\cite{ovsjanikov2012functional}) arise in this application as they allow to align features in the spectral domain and to extract a point to point map as a post processing step. With respect to CDMD, Eynard et al.~\cite{eynard2016coupled} investigate a close variant of our CDMD problem, and solved it directly using a non-linear conjugate gradients approach. An alternative formulation was studied in \cite{huang2014functional}, based on the observation that the matrix 
\[
Z = \begin{pmatrix} I & A \\ B & I \end{pmatrix}
\]
is low-rank when $AB = I$. Instead of minimizing the rank of $Z$, Huang et al.~\cite{huang2014functional} replace the low-rank constraint with its convex relaxation expressed via the nuclear norm~\cite{candes2011robust}.

Our approach depends on the following straightforward insight. Under the change of variables $B$, the energy functional in~\cref{eq:fb_pbm2} becomes fully \emph{separable}. Namely, if we denote 
\begin{align}
\hat{f}(A) = \frac{1}{2} |AX-Y|_F^2 \ , \quad \tilde{f}(B) = \frac{1}{2} |X-BY|_F^2 \ ,    
\end{align}
then we seek to minimize $\hat{f}(A) + \tilde{f}(B)$ subject to the constitutive invertibility constraints. This understanding calls for the development of an Alternating Direction Method of Multipliers (ADMM)-type approach~\cite{boyd2011distributed}. ADMM is advantageous in effectively solving separable optimization problems, since it systematically leads to splitting schemes composed of potentially simpler minimization tasks. Moreover, the theory associated with ADMM-based techniques is well-developed with several general results related to convergence, optimality conditions and stopping criteria. Unfortunately, the constraints associated with our problem are \emph{non-linear}, and thus while one can employ an ADMM approach, the theoretical guarantees of standard ADMM do not apply. Recently, Gao et al.~\cite{gao2018admm} showed that under mild assumptions, ADMM with \emph{multiaffine} constraints converges if the penalty parameter in the augmented Lagrangian is sufficiently large. In~\Cref{sec:converging_cdmd}, we show that CDMD can be modified to fit a family of optimization problems that are considered in~\cite{gao2018admm} for which converging ADMM schemes can be devised.

\subsection{A splitting scheme} 

We now turn to present the main algorithm in this work. Our starting point is to define the augmented Lagrangian for problem~\cref{eq:fb_pbm2} given by,
\begin{align} \label{eq:aug_lag}
\mca{L}(A,B,Q) = \hat{f}(A) + \tilde{f}(B) + \frac{\rho}{2} |R(A,B) + Q|_F^2 - \frac{\rho}{2} |Q|_F^2 \ ,
\end{align}
where $\rho \in \mbb{R}^{+}$ is a scalar penalty parameter, the matrix $R(A,B)$ combines the constitutive constraints into a single matrix, and the matrix $Q$ is the scaled dual variable (see e.g.,~\cite[Section 3.1.1]{boyd2011distributed}). Specifically, the matrices $R$ and $Q$ are given by
\[
R(A,B) = \begin{pmatrix} AB-I \\ BA-I \end{pmatrix} \in \mbb{R}^{2r\times r} \ , \quad Q = \begin{pmatrix} Q_1 \\ Q_2 \end{pmatrix} \in \mbb{R}^{2r\times r} \ .
\]

We note that if one adopts the method of multipliers approach, the augmented Lagrangian $\mca{L}(A,B,Q)$ could be directly minimized, as was done in~\cite{eynard2016coupled}. However, the term $|R(A,B)+Q|_F^2$ includes a quartic combination of unknowns, and thus the optimization problem~\cref{eq:aug_lag} is highly non-linear. Instead, our numerical scheme splits the updates so that $A$ and $B$ are not updated jointly but in an alternate fashion. Specifically, given initial $A^0,B^0,Q^0$ and $\rho$, ADMM takes the form of
\begin{enumerate}
    \item $ A^{k+1} = \argmin_A \hat{f}(A) + \frac{\rho}{2} |R(A,B^k)+Q^k|_F^2$
    \vspace{3mm}
    \item $ B^{k+1} = \argmin_B \tilde{f}(B) + \frac{\rho}{2} |R(A^{k+1},B)+Q^k|_F^2$
    \vspace{3mm}
    \item $ Q^{k+1} = Q^k + R(A^{k+1},B^{k+1}) $
\end{enumerate}
Below, we show that minimizing Steps~($1$) and ($2$) lead in both cases to a Sylvester Equation which can be efficiently solved using the QR decomposition, see~\cite{bartels1972solution} for further details. The update in Step~($3$) is trivial and requires a single evaluation of $R$. Overall, we obtain an efficient algorithm with time complexity of $\mca{O}(K r^3)$, where $K$ is the total number of iterations.

The minimization tasks in Steps~($1$) and ($2$) are relatively simple as they comprise of energy functionals that are quadratic in $A$ and in $B$, respectively. Thus, the associated first order optimality conditions are $\emph{linear}$. For instance, the Jacobian of the energy in Step~($1$) is
\begin{align*}
\nabla_A \left[ \mca{L}(A,B^k,Q^k) \right] &= \nabla_A \, \hat{f}(A) + \frac{\rho}{2} \nabla_A \left( | R(A,B^k) + Q^k |_F^2 \right)  \\
&= (AX-Y)X^T + \rho \left(A\,B^k - I + Q_1^k \right) (B^k)^T + \rho (B^k)^T \left(B^k A - I + Q_2^k \right) \ .
\end{align*}
After re-arrangement and equating to zero, we arrive at the following \emph{Sylvester Equation}, $C_1 A + A \, C_2 = C_3$, which is linear in $A$. The matrices $C_1, C_2$ and $C_3$ are given by
\begin{align} \label{eq:upd_sylmat_a} \begin{split}
    C_1 &= \rho (B^k)^T B^k \ , \\
    C_2 &= XX^T + \rho B^k (B^k)^T \ , \\
    C_3 &= YX^T + 2 \rho (B^k)^T - \rho Q_1^k (B^k)^T - \rho (B^k)^T Q_2^k \ .
\end{split} \end{align}
The derivation for Step~($2$) follows along the same lines, yielding a different Sylvester Equation $D_1 B + B D_2 = D_3$ with coefficient matrices given by
\begin{align} \label{eq:upd_sylmat_b} \begin{split}
    D_1 &= \rho (A^{k+1})^T A^{k+1} \ , \\
    D_2 &= YY^T + \rho A^{k+1} (A^{k+1})^T \ , \\
    D_3 &= XY^T + 2 \rho (A^{k+1})^T - \rho (A^{k+1})^T Q_1^k - \rho Q_2^k (A^{k+1})^T \ .
\end{split} \end{align}

\subsection{The numerical algorithm} 

We summarize our technique for computing consistent dynamic mode decomposition in~\cref{alg:cdmd}. Note that Steps~$1-2$ and $10-11$ are shared with~\cref{alg:exact_dmd}, whereas our main contribution is provided in Steps~$3-9$ where the construction of the DMD matrix $A$ is described. We note that the algorithm below describes how to compute an approximation of the forward dynamics $A$ and its associated decomposition, however, an estimate of the backward dynamics can be extracted as well by defining $B = B^k$, where $k$ is the last iteration index.

\begin{algorithm}
\caption{Consistent Dynamic Mode Decomposition (\code{CDMD})}
\label{alg:cdmd}
\begin{algorithmic}[1]
\STATE{Input matrices $\tilde{X},\tilde{Y} \in \mbb{R}^{m\times n}$ and scalars $r, \rho \in \mbb{R}$}
\vspace{2mm}
\STATE{Compute the SVD of $\tilde{X}=\tilde{U} \tilde{S} \tilde{V}^*$, and generate $X = \tilde{U}_r^* \tilde{X}, Y = \tilde{U}_r^* \tilde{Y}$}
\vspace{2mm}
\STATE{Initialize $A^0 = YX^{+} , B^0 = XY^{+}, Q^0 = 0$}
\FOR{$k=0,1,2,...$}
\STATE{Solve $A^{k+1} = \code{sylvester}(C_1, C_2, C_3)$, using Eq.~\cref{eq:upd_sylmat_a} }
\vspace{1mm}
\STATE{Solve $B^{k+1} = \code{sylvester}(D_1, D_2, D_3)$, using Eq.~\cref{eq:upd_sylmat_b} }
\vspace{1mm}
\STATE{Update $Q^{k+1} = Q^k + R(A^{k+1}, B^{k+1})$}
\STATE{Update $\rho$ following Eq.~\cref{eq:update_rho}}
\ENDFOR
\vspace{2mm}
\STATE{Compute the EIG of the last $A$, with $A v_j = \lambda_j v_j$, where $v_j \in \mbb{C}^r, \; \lambda_j \in \mbb{C}$}
\vspace{2mm}
\STATE{The DMD spectrum is defined as the set of eigenvalues $\lambda_j$, and vectors $\psi_j = \lambda_j^{-1} \tilde{Y} \tilde{V}_r \tilde{S}_r^{-1} v_j$}
\end{algorithmic}
\end{algorithm}

\subsection{Stopping criteria}

To establish a practical stopping condition, we keep track of two residual quantities that are related to the primal and dual problems. A similar termination approach is described in~\cite{boyd2011distributed}. We define the following \emph{primal residual} and \emph{dual residual},
\begin{align}
    r^k = R(A^k,B^k) \ , \quad s^k = \rho \begin{pmatrix} A^{k}-A^{k-1} \\ B^{k}-B^{k-1} \end{pmatrix} \ ,
\end{align}
where the termination rule we employ is given by $|r^k|_F \leq \epsilon^{\text{pri}}$ and $|s^k|_F \leq \epsilon^{\text{dual}}$. The tolerances $\epsilon^{\text{pri}}$ and $\epsilon^{\text{dual}}$ can be computed using absolute and relative thresholds, such as
\begin{align*} \begin{split}
    \epsilon^{\text{pri}} &= \sqrt{r} \epsilon^{\text{abs}} + \epsilon^{\text{rel}} \max \left\{ |A^k B^k |_F, |B^k A^k |_F \right\} \ , \\
    \epsilon^{\text{dual}} &= \sqrt{2r} \epsilon^{\text{abs}} + \epsilon^{\text{rel}} \rho|Q^k|_F \ .
\end{split} \end{align*}
                    
\subsection{Dynamic update of the penalty parameter $\rho$}
\label{subsec:vary_rho}

In general, varying $\rho$ based on the current estimates of the primal and dual residuals may lead to faster convergence rates. We implement a simple scheme that was proposed in e.g.,~\cite{boyd2011distributed} and is given by
\begin{align} \label{eq:update_rho}
    \rho^{k+1} := \left\{ 
                \begin{array}{ll}
                  \tau \rho^k \quad\quad & \text{if } |r^k|_F > \mu |s^k|_F \\
                  \rho^k/\tau  & \text{if } |s^k|_F > \mu |r^k|_F \\
                  \rho^k  & \text{otherwise},
                \end{array} 
                \right.
\end{align}
where we take $\tau = 2$ and $\mu = 5$ in practice.




\begin{figure*}[h]
 \centering
 \includegraphics[width=\linewidth]{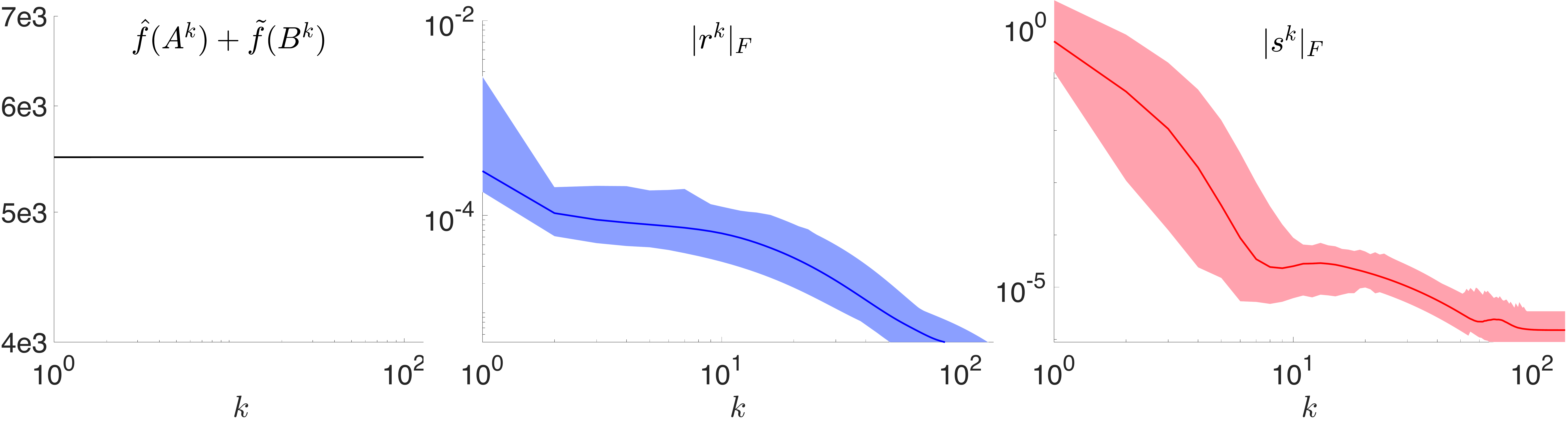}
 \vspace{-8mm}
 \caption{The empirical convergence of our \Cref{alg:cdmd} when applied to non-linear data with high levels of noise is demonstrated in the above plot. Our method terminates in $\approx 100$ steps, where the objective function is stabilized on an optimal value (left) and the primal and dual residuals converge rapidly (middle and right). We repeat this test $N=1000$ times, and we show the variance in convergence via the shaded areas where the average is represented by bold curves.}
 \label{fig:cdmd_stats}
\end{figure*}

\section{Provably Convergent CDMD Scheme}
\label{sec:converging_cdmd}

Unfortunately, while the above~\cref{alg:cdmd} is effective and behaves well in practice as we show in~\Cref{sec:results} and in Fig.~\ref{fig:cdmd_stats}, it is not \emph{provably} convergent. In what follows, we address this shortcoming and propose an alternative converging scheme, which requires only an additional negligible amount of computations. To this end, we follow the recent work of Gao et al.~\cite{gao2018admm} which showed that under certain conditions, ADMM and its convergence can be extended to include multiaffine constraints. In particular, we show that by introducing additional variables to the CDMD problem~\cref{eq:fb_pbm2}, the obtained minimization problem is of the required form, while satisfying all the necessary conditions in~\cite{gao2018admm}.

Gao et al. investigate the convergence of ADMM for problems taking the form,
\begin{align} \label{eq:gao_admm}
\minz_{\mca{A,B,C}} \quad h(\mca{A,B,C}) \ , \quad \text{s.t.} \quad \mca{P}(\mca{A,B}) + \mca{Q}(\mca{C}) = 0 \ ,
\end{align}
where $\mca{A} = (A_0,A_1,...,A_{n_a})$, $\mca{B} = (B_0,B_1,...,B_{n_b})$, and a variable block $\mca{C}$. In addition, we have that $h(\mca{A,B,C}) = f(\mca{A,B})+g(\mca{C})$. Finally, $\mca{Q}$ is a linear map and, in contrast to ``standard'' ADMM problems, $\mca{P}$ is a \emph{multiaffine} map. Namely, the transformation obtained from fixing all variables $A_i$ and $B_j$ but one, is affine. It is shown in~\cite{gao2018admm} that when several assumptions on $h,\mca{P,Q}$ are met, an ADMM scheme converges to a constrained stationary point, i.e., the sequence $\{\mca{A}^k,\mca{B}^k,\mca{C}^k\}_{k=0}^\infty$ is bounded, and that every limit point $(\mca{A}^*,\mca{B}^*,\mca{C}^*)$ is a constrained stationary point. While various configurations of assumptions are considered in~\cite{gao2018admm}, we list here a more restrictive set of conditions that hold in our case.
\begin{assumption*} \label{assum:admm} Solving problem~\cref{eq:gao_admm}, the following hold.
\begin{enumerate}
    \item The update order is $A_0,A_1,...,A_{n_a},B_0,B_1,...,B_{n_b}$ and a single block $\mca{C}$.
    
    \item $\im(\mca{Q}) \supseteq \im(\mca{P})$.
    
    \item The objective $h(\mca{A,B,C})$ is coercive on the feasible set $$ \Omega = \left\{(\mca{A,B,C}) : \mca{P}(\mca{A,B})+\mca{Q}(\mca{C}) = 0 \right\} \ .$$
    
    \item The function $f(\mca{A,B})$ can be written as $$ f(\mca{A,B}) = \sum_i^{n_a} \hat{f}_i(A_i) + \sum_j^{n_b} \tilde{f}_j(B_j) \ ,$$ where every $\hat{f}_i$ and $\tilde{f}_j$ are $(m_i,M_i)$- and $(m_j,M_j)$-strongly convex functions.
    
    \item The function $g(\mca{C})$ is a $(m,M)$-strongly convex function.
    
    \item For sufficiently large penalty $\rho$, every ADMM subproblem attains its optimal value.
    
\end{enumerate}
\end{assumption*} 

To motivate our discussion, we present an illustrative example related to Nonnegative Matrix Factorization (NMF). As we show below, this problem is similar to ours with respect to the biaffine constraints, and thus it provides a natural starting point for our case. Given a matrix $Z$, its NMF involves the task of finding a pair of nonnegative matrices $A \geq 0$ and $B \geq 0$ such that $Z = AB$~\cite{lee2001algorithms}. An ADMM formulation to NMF was originally proposed in~\cite{boyd2011distributed}, yielding the following problem,
\begin{align} \label{eq:nmf}
\minz_{A,B,C} \quad \imath(A) + \imath(B) + \frac{1}{2}|C-Z|_F^2 \ , \quad \text{s.t.} \quad C = AB \ , 
\end{align}
    where $\imath$ is the indicator function, i.e., $\imath(A)=0$ if $A \geq 0$ and $\imath(A) = \infty$ otherwise. Gao and colleagues reformulate~\cref{eq:nmf} to arrive at an optimization problem whose subproblems are easy to solve while meeting the assumptions required for convergence. The modified version is given by
\begin{align} \begin{split} \label{eq:nmf2}
\minz_{A,A',B,B',C,A'',B''} & \quad \imath(A') + \imath(B') + \frac{1}{2}|C-Z|_F^2 + \frac{\mu}{2}|A''|_F^2 + \frac{\mu}{2}|B''|_F^2 \ , \\
\text{subject to} & \quad C = AB, \; A = A'+A'', \; B = B'+B''  \ . 
\end{split} \end{align}
The update order of the variables is $B,B',A,A'$ and $(C,A'',B'')$. We stress that problem~\cref{eq:nmf2} satisfies a different set of assumptions than those appear in~\cref{assum:admm}, but it is well within the family of problems considered in~\cite{gao2018admm}. We refer to their paper for additional details of the NMF problem considered in relation to converging ADMM schemes.

We now turn to modify the CDMD problem~\cref{eq:fb_pbm2} to a form which fits all the conditions in~\cref{assum:admm} and thus its ADMM is provably convergent, due to~\cite{gao2018admm}. We observe that our invertibility constraints $AB=I$ and $BA=I$ are reminiscent of the NMF constraints, and, in particular, they are biaffine with respect to $(A,B)$. Moreover, our objective function consists of highly smooth Frobenius norm terms. Encouraged by these similarities, we introduce the auxiliary variables $C,A',A'',B',B''$, and we modify the above~\cref{eq:fb_pbm2} to arrive at the following minimization,
\begin{align} \begin{split}
\minz_{A,A',B,B',C,A'',B''} & \quad \frac{1}{2} |A'X-Y|_F^2 + \frac{1}{2} |X-B'Y|_F^2 + \frac{\nu}{2}|C-I|_F^2 + \frac{\mu}{2} |A''|_F^2 + \frac{\mu}{2} |B''|_F^2 \ , \\
\text{subject to} & \quad C = AB, \, C = BA, \, A = A'+A'', \, B = B'+B'' \ ,
\end{split} \label{eq:fb_pbm3} \end{align}
where $\nu,\mu \in \mbb{R}^+$ are penalty parameters for the $C, A''$ and $B''$ variables.

\begin{algorithm}
\caption{Provably Convergent CDMD (\code{CDMD2})}
\label{alg:cdmd2}
\begin{algorithmic}[1]
\STATE{Input matrices $\tilde{X},\tilde{Y} \in \mbb{R}^{m\times n}$ and scalars $r, \rho, \mu \in \mbb{R}$}
\vspace{2mm}
\STATE{Compute the SVD of $\tilde{X}=\tilde{U} \tilde{S} \tilde{V}^*$, and generate $X = \tilde{U}_r^* \tilde{X}, Y = \tilde{U}_r^* \tilde{Y}$}
\vspace{2mm}
\STATE{Initialize $A^0 = A'^0 = YX^{+} , B^0 = B'^0 = XY^{+}, A''^0 = B''^0 = 0, C^0 = I, Q^0 = 0$}
\FOR{$k=0,1,2,...$}
\STATE{Solve $A^{k+1} = \code{sylvester}(A_1, A_2, A_3)$, where 
\begin{align*}
A_1 &= I+(B^k)^T B^k \ , \\
A_2 &= B^k (B^k)^T \ , \\ 
A_3 &= (C^k-Q_1^k)(B^k)^T + (B^k)^T(C^k-Q_2^k) + A'^k + A''^k - Q_3^k \ .
\end{align*} }
\vspace{-3mm}
\STATE{Solve $A'^{k+1} = \code{linsolve}\left(\rho I + XX^T, YX^T + \rho(A^{k+1}-A''^k+Q_3^k) \right)$ }
\vspace{1mm}
\STATE{Solve $B^{k+1} = \code{sylvester}(B_1, B_2, B_3)$, where 
\begin{align*}
B_1 &= I + (A^{k+1})^T A^{k+1} \ , \\ 
B_2 &= A^{k+1} (A^{k+1})^T \ , \\
B_3 &= (A^{k+1})^T(C^k-Q_1^k) + (C^k-Q_2^k)(A^{k+1})^T + B'^k + B''^k - Q_4^k \ .
\end{align*} }
\vspace{-3mm}
\STATE{Solve $B'^{k+1} = \code{linsolve}\left(\rho I + YY^T, XY^T + \rho(B^{k+1}-B''^k+Q_4^k) \right)$ }
\vspace{1mm}
\STATE{Solve $C^{k+1} = \frac{\rho}{2\rho+\nu}(A^{k+1}B^{k+1}+B^{k+1}A^{k+1} +Q_1^k+Q_2^k) + \frac{\nu}{2\rho+\nu} I$ }
\vspace{1mm}
\STATE{Solve $A''^{k+1} = \frac{\rho}{\mu+\rho} (A^{k+1}-A'^{k+1} + Q_3^k)$ }
\vspace{1mm}
\STATE{Solve $B''^{k+1} = \frac{\rho}{\mu+\rho} (B^{k+1}-B'^{k+1} + Q_4^k)$ }
\vspace{1mm}
\STATE{Update $Q^{k+1} = Q^k + \mca{R}(\mca{A}^{k+1}, \mca{B}^{k+1},\mca{C}^{k+1})$}
\STATE{Update $\rho$ following Eq.~\cref{eq:update_rho} }
\ENDFOR
\vspace{2mm}
\STATE{Execute steps $(10)-(11)$ of~\Cref{alg:cdmd}}
\end{algorithmic}
\end{algorithm}

To verify that~\cref{eq:fb_pbm3} meets all the required conditions, we denote $\hat{f}(A') = \frac{1}{2} |A'X-Y|_F^2$, $\tilde{f}(B') = \frac{1}{2} |X-B'Y|_F^2$, and $g(C,A'',B'')=\frac{\nu}{2}|C-I|_F^2 + \frac{\mu}{2} |A''|_F^2 + \frac{\mu}{2} |B''|_F^2$. Also, we define the following residual
\begin{align*}
\mca{R}(\mca{A,B,C}) = \mca{P}(A,A',B,B') + \mca{Q}(C, A'', B'') = \begin{pmatrix} AB \\ BA \\ A-A' \\ B-B' \end{pmatrix} + \begin{pmatrix} -C \\ -C \\ -A'' \\ -B'' \end{pmatrix} \ .
\end{align*}
The conditions in~\cref{assum:admm} hold because the update order is $A,A',B,B'$ and $(C,A'',B'')$ as we show below in~\cref{alg:cdmd2}. The image of $\mca{Q}$ is indeed a superset of $\mca{P}$'s image, since it is the (minus) identity transformation in each of its entries, and thus span the entire space. The objective function $h$ is coercive on the feasible set, because its terms behave as $|x|_F^2$, and therefore whenever $|x|_F \rightarrow \infty$ so does $|x|_F^2$. Under some mild conditions, namely, that $X$ and $Y$ are full rank matrices, the function $f$ is composed of $(m,M)$-strongly convex functions as we show in~\cref{app:f_sconvex}. Similarly, $g$ is a strongly convex function because the Hessian of its terms is positive definite. Finally, the subproblems in our formulation are trivial, linear or a Sylvester-type equation and thus attain their optimal value when $\rho$ is sufficiently large.

We conclude this section with presenting our convergent ADMM scheme along with the specification of its subproblems. The derivation of the matrix expressions that take part in lines $5$ and $7$ could be carried over in a fashion similar to Eqs.~\cref{eq:upd_sylmat_a} and \cref{eq:upd_sylmat_b}. We note that lines $6$ and $8$ of \Cref{alg:cdmd2} involve a call to $X=\code{linsolve}(A,B)$ which numerically solves the system $XA=B$. 

\section{Results}
\label{sec:results}

In this section, we evaluate the proposed \code{CDMD} approach and compare it to several state-of-the-art techniques for computing DMD matrices. In particular, we compare against \code{Exact DMD}~\cite{tu2013dynamic}, \code{fbDMD}~\cite{dawson2016characterizing}, \code{tlsDMD}~\cite{hemati2017biasing} and \code{optimized DMD}~\cite{askham2018variable}. The dynamical systems we consider appeared previously e.g., in~\cite{dawson2016characterizing,askham2018variable}, and thus can be considered as ``benchmark'' examples for quantitative and qualitative study of DMD algorithms.

\begin{figure*}[h]
 \centering
 \includegraphics[width=\linewidth]{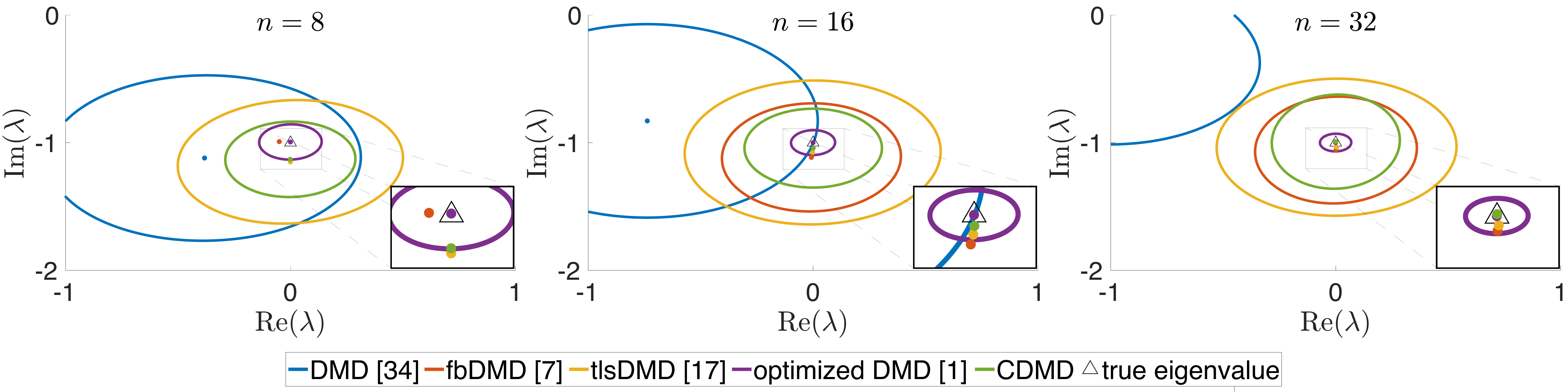}
 \vspace{-3mm}
 \caption{We plot $95\%$ confidence ellipses (see~\cite{dawson2016characterizing}) for estimating one of the eigenvalues of a periodic linear system~\cref{eq:linp_nn} when varying number of observations $n=8,16,32$ are given. The zoom in boxes show the average estimation for each method. The results above indicate that \code{CDMD} is second to best in terms of accuracy and variance for all values of $n$.}
 \label{fig:linp_nn_eigvals}
 \vspace{-3mm}
\end{figure*}

\subsection{A periodic linear system}

In this example, we use the following linear and non-normal system
\begin{align} \label{eq:linp_nn}
    \dot{z} = \begin{pmatrix} 1 & -2 \\ 1 & -1 \end{pmatrix} z \ ,
\end{align}
where the system has purely imaginary eigenvalues that are given by $\lambda = \pm i$. Eq.~\cref{eq:linp_nn} is integrated over the $[0,2\pi]$ temporal segment, starting from the initial point $z_0 = [1 \; 0.1]^T$. To stress test our method, we investigate this system when relatively low number of observations is given and high levels of white Gaussian noise affect the data. Specifically, we show in~\cref{fig:linp_nn_eigvals} the performance of various methods for computing the eigenvalue $-i$ when noise with variance $\sigma^2 = 0.1$ and Signal-to-Noise (SNR) ratio of $8.6\ \mathrm{dB}$ is introduced. We repeat our experiment $N=10^4$ times, and the average of each of the methods is marked by a dot with a corresponding color. Additionally, we plot the ellipses which enclose the region of $95\%$ of the estimates that are closest to the true eigenvalue for each of the techniques. We use the values $n=8,16,32$ for the number of observations, which make the system overdetermined as it is two-dimensional. Nevertheless, these values are relatively small in comparison to related work on this example, see e.g.,~\cite{dawson2016characterizing}. 

\begin{figure*}[t]
 \centering
 \includegraphics[width=\linewidth]{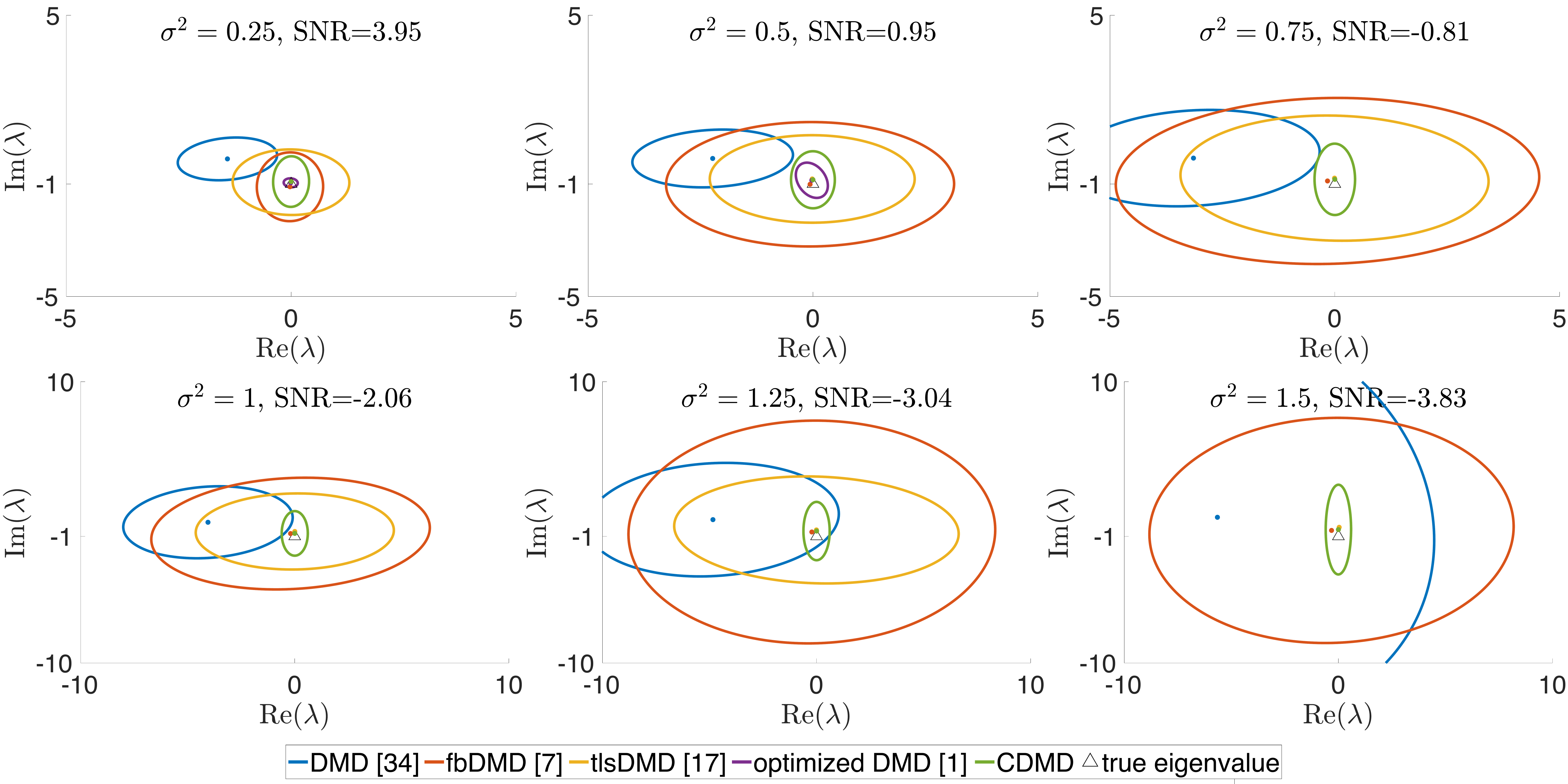}
 \vspace{-3mm}
 \caption{We compare the above methods based on their approximation for the eigenvalue of System~\cref{eq:linp_nn} when various levels of noise are introduced, $-4 \leq \text{SNR} \leq 4$. Interestingly, while \code{optimized DMD} is extremely accurate when $\text{SNR}>0$, it fails for higher levels of noise, and thus it does not appear on these graphs. The methods \code{fbDMD} and \code{tlsDMD} perform well in terms of average, but their spread is much larger than our results which maintain relatively small spread as well as accurate average.}
 \label{fig:linp_nn_eigvals_snr}
\end{figure*}

Overall, \code{optimized DMD} achieves excellent results in terms of spread and average values, across all values of $n$. On the other end, \code{exact DMD} struggles both in accuracy and spread. \code{fbDMD} and \code{tlsDMD} exhibit comparable performance, except for $n=8$ where \code{fbDMD} produces a correct mean, but with an extremely large deviation. Finally, our approach outputs consistent deviation and averages, regardless of the value of $n$. We additionally experiment with various high level of noise $-4 \leq \text{SNR} \leq 4$ and present the results in~\cref{fig:linp_nn_eigvals_snr}. Note that the bottom row axes are twice as large as the axes in the top row. As can be seen in the graphs, \code{optimized DMD} is very accurate as long as $\text{SNR} > 0$, but fails when the signal-to-noise ratio drops below zero, and therefore it is omitted from the other graphs. In most cases, \code{Exact DMD} produces poor approximations when compared to the other methods. In comparison, \code{fbDMD} and \code{tlsDMD} generate estimates that are centered around the eigenvalue in general, with growing spread as the SNR decreases. Remarkably, our approach exhibits the least increase in deviation when compared to all other techniques, while producing a relatively accurate average.

\begin{figure*}[t]
 \centering
 \includegraphics[width=\linewidth]{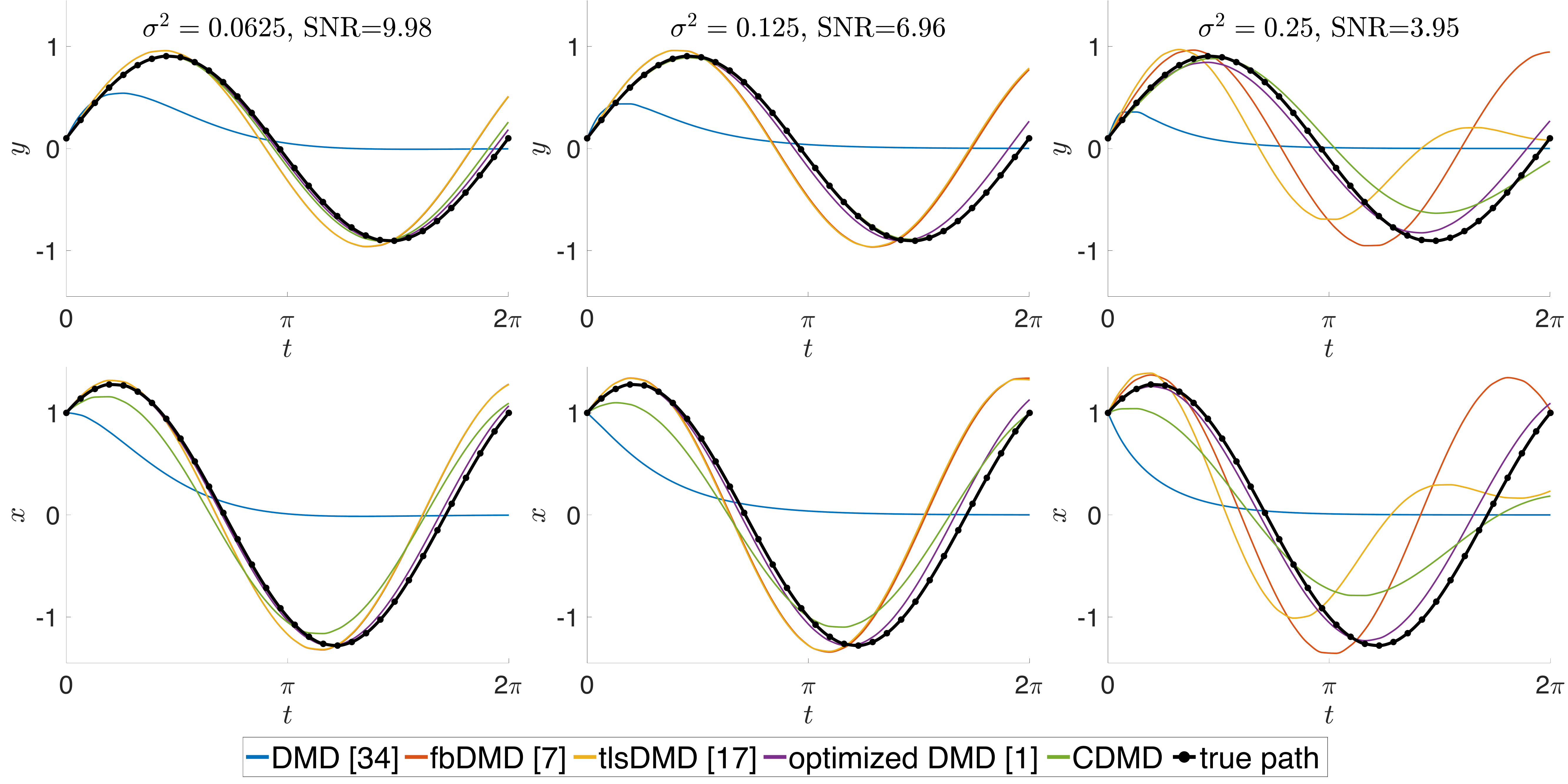}
 \vspace{-8mm}
 \caption{We reconstruct the trajectory of the periodic system~\eqref{eq:linp_nn} using the computed DMD matrices for various noise variances with $32$ observations. Most methods yield paths that are close to the true trajectory, where \code{optimized DMD} and our method obtain the best results.}
 \label{fig:linp_nn_op}
\end{figure*}

In addition, we reconstruct the trajectory using the approximations of the dynamics provided by each of the methods, and we plot the results in Fig.~\ref{fig:linp_nn_op} separated to $y$-coordinate (top row) and $x$-coordinate (bottom row) over time. It is evident that \code{Exact DMD} yields a highly distorted path, whereas the other methods are generally close to the true trajectory. As the amount of noise increases, \code{fbDMD} and \code{tlsDMD} develop a significant shift in phase. We measure the distance between the computed paths to the desired curve and we observe that our method achieves second to best results after \code{optimized DMD}. Specifically, for $\sigma^2 = 0.125$, the $L_2$ error between the computed path to the ground-truth trajectory divided by the length of the latter is $0.0837$ and $0.2611$ for \code{optimized DMD} and \code{CDMD}, respectively. When $\sigma^2 = 0.25$, the error is $0.1403$ and $0.7844$ for \code{optimized DMD} and \code{CDMD}. In comparison, the other methods yield errors that are five times larger or more.

\subsection{Dominant and hidden dynamics}

\begin{figure*}[t]
 \centering
 \includegraphics[width=\linewidth]{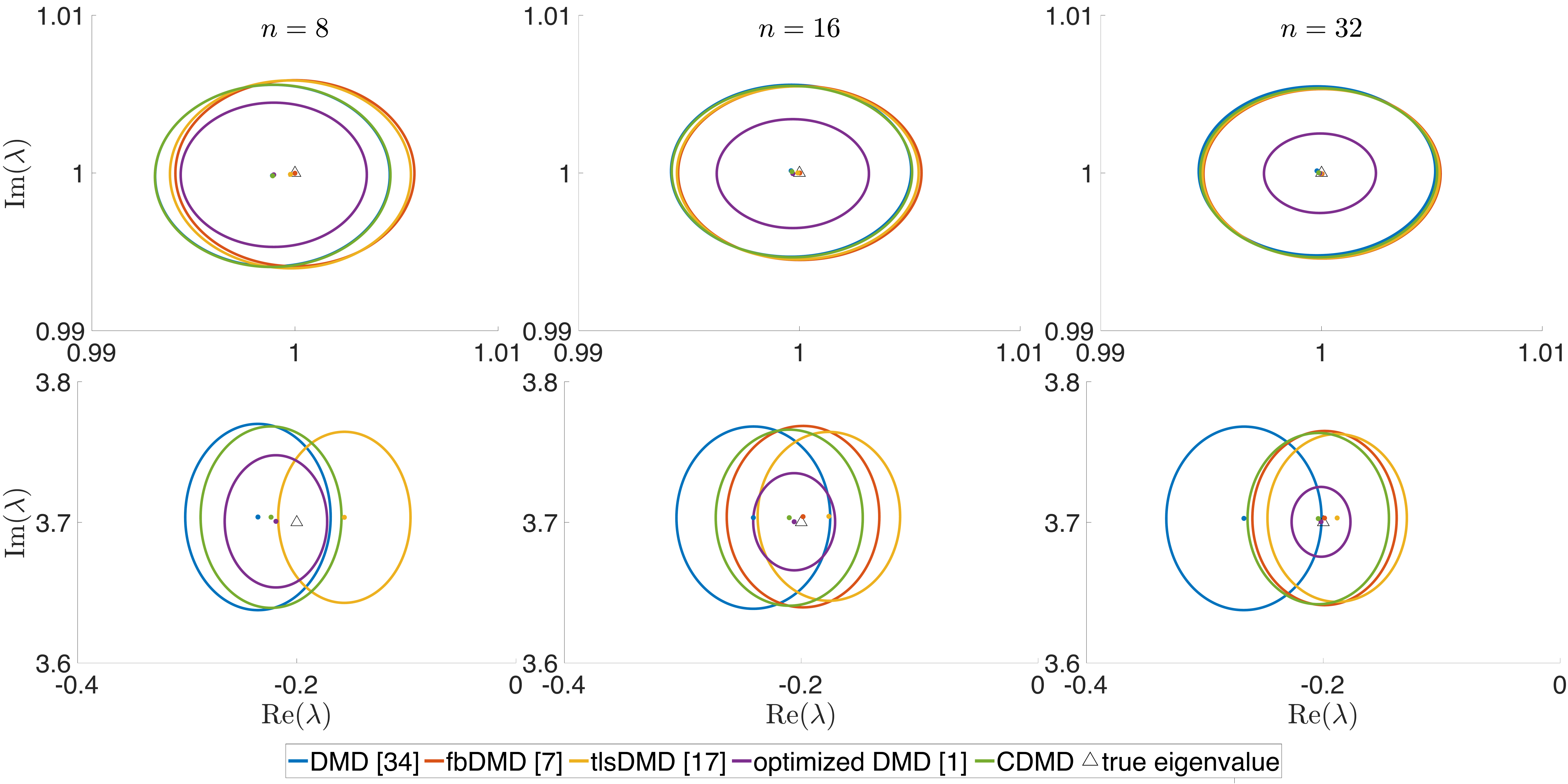}
 \vspace{-7mm}
 \caption{Given a noisy superposition of sine functions, we estimate the system's eigenvalues for various number of observations $n$ and noise with variance $\sigma^2=0.25$ and $\text{SNR}=30\ \mathrm{dB}$. We observe that the dominant eigenvalue is approximated well (top row), whereas the hidden dynamics is achieved by most methods with decreasing error as $n$ grows (bottom row).}
 \label{fig:sine_eigvals}
\end{figure*}

The next system is a superposition of a growing sine function and a decaying sine function given by
\begin{align} \label{eq:nonlin_sine}
    z(x,t) = \sin(k_1 x-\omega_1 t)\exp(\gamma_1 t) + \sin(k_2 x-\omega_2 t)\exp(\gamma_2 t) \ ,
\end{align}
where in our experiments we used $k_1=1, \, \omega_1=1, \, \gamma_1=1$ and $k_2=0.4, \, \omega_2=3.7, \, \gamma_2=-0.2$. This example is more challenging than the previous one since it involves dynamical features which are of lower magnitude alongside dominant structures. The eigenvalues of this system are of the form $\gamma_i \pm \omega_i, i = 1,2$, where the ``dominant'' mode is associated with $i=1$ and the ``hidden'' mode is linked to $i=2$. In~\Cref{fig:sine_eigvals}, we compute $N=10^4$ times the eigenvalues of the system while employing a noise level of $\sigma^2=0.25$, $\text{SNR}=30\ \text{dB}$ over the observations. The results show that for the dominant dynamics, most methods perform well where \code{optimized DMD} obtains improved estimates as $n$ increases (top row). For the hidden mode, similar results are obtained for $n= 16, 32$, whereas for the lowest $n=8$, \code{fbDMD} does not appear in the plot and \code{tlsDMD} is shifted differently than the other approaches (bottom row).

\begin{figure*}[h]
 \centering
 \includegraphics[width=\linewidth]{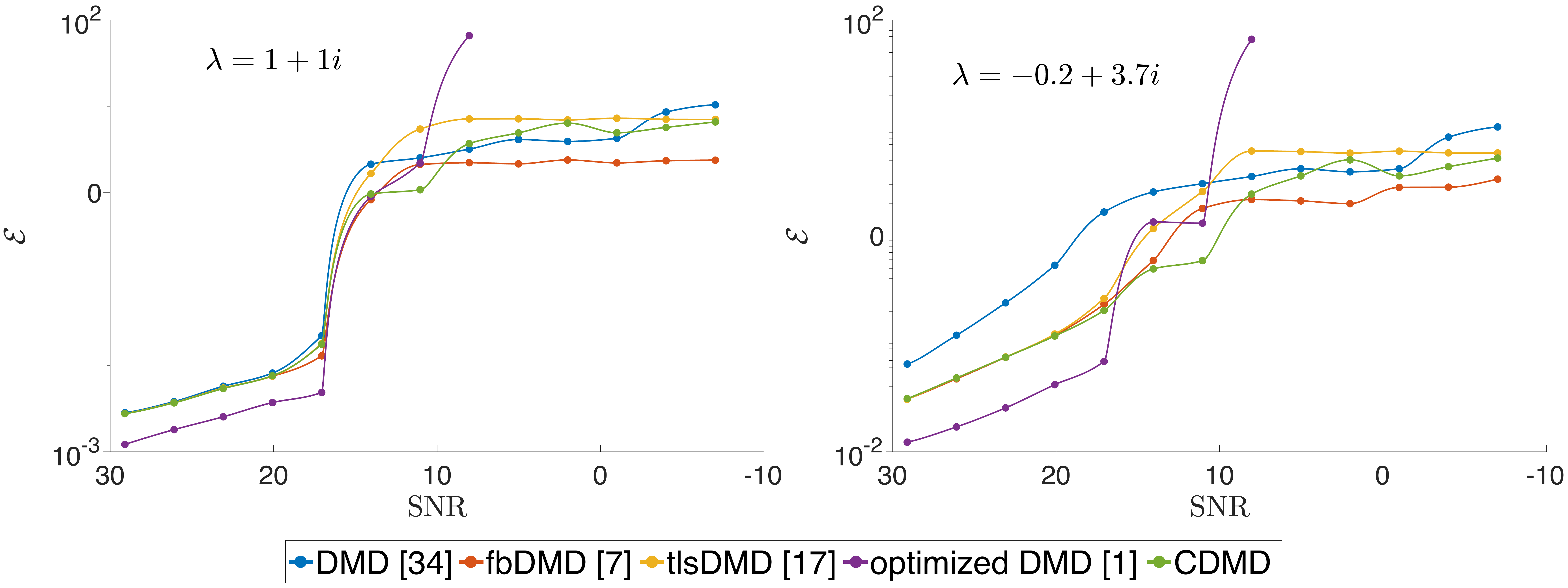}
 \vspace{-7mm}
 \caption{We compute approximations of the dominant and hidden DMD eigenvalues for various levels of noise, $-10 \leq \text{SNR} \leq 30$ and we compute the error for each method using Eq.~\eqref{eq:err_metric}. Naturally, for low levels of noise, most methods perform reasonably well, whereas when $\text{SNR} \leq 17$ the spread becomes orders of magnitude larger. See the text for further details.}
 \label{fig:sine_eigvals_snr}
\end{figure*}

In addition, we investigate this system across different levels of noise. In particular, we set $\sigma^2 = 2^{-2},\ 2^{-1},..., 2^{10}$ corresponding to SNR in the range $[-10,30]$. Each noise level is used $N=10^3$ times, for which we compute both the dominant and hidden DMD eigenvalues. We show the error results of the different methods in Fig.~\cref{fig:sine_eigvals_snr}, where the error is a linear combination of the average error between the computed eigenvalue and the ground-truth and the minimum radius of the deviation ellipse. Formally, 
\begin{align} \label{eq:err_metric}
\mathcal{E} = a |\lambda_{\text{avg}}-\lambda_{\text{gt}}| + (1-a) r_{min} \ ,
\end{align}
where $\lambda_{\text{avg}}$ is the average taken over all eigenvalue estimates, $\lambda_{\text{gt}}$ is the analytic eigenvalue, and $r_{min}$ is the minimum radius. In our experiments, we used $a=0.9$. Similar to Fig.~\cref{fig:linp_nn_eigvals_snr}, when SNR approaches zero, \code{optimized DMD} fails and thus its graphs are shorter. Interestingly, up to a certain SNR, all methods present similar error behavior, where at SNR $\approx 17$ there is an exponential increase in the error estimates. When inspecting the individual results, it seems like this high level of noise leads to an extremely large deviation in results, which further affects our error measure.

\subsection{Cylinder wake}

\begin{figure*}[h]
 \centering
 \includegraphics[width=\linewidth]{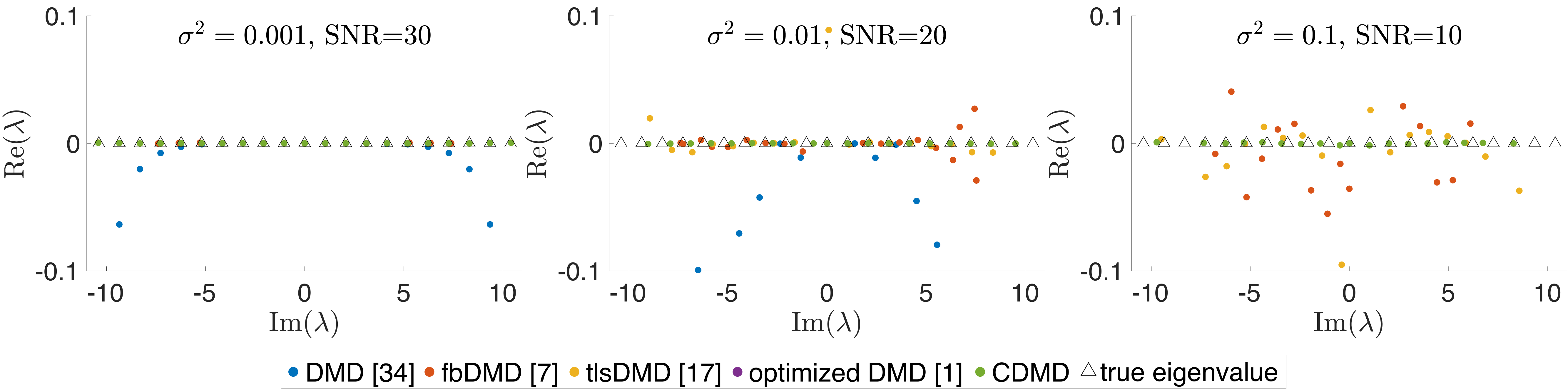}
 \vspace{-8mm}
 \caption{We compute $r=21$ eigenvalues using simulation snapshots of a cylinder wake for noise levels $\sigma^2=0.001,0.01,0.1$ respectively $\text{SNR}= 30, 20, 10\ \text{dB}$. As the noise increases, our method maintains its zero growth estimate (notice that the $x$-axis represents the real part, cf. Fig.~\cref{fig:sine_eigvals}), whereas the other methods produce significant erroneous growth/decay estimates.}
 \label{fig:cyn_eigvals}
\end{figure*}

\begin{figure*}[t]
 \centering
 \includegraphics[width=.5\linewidth]{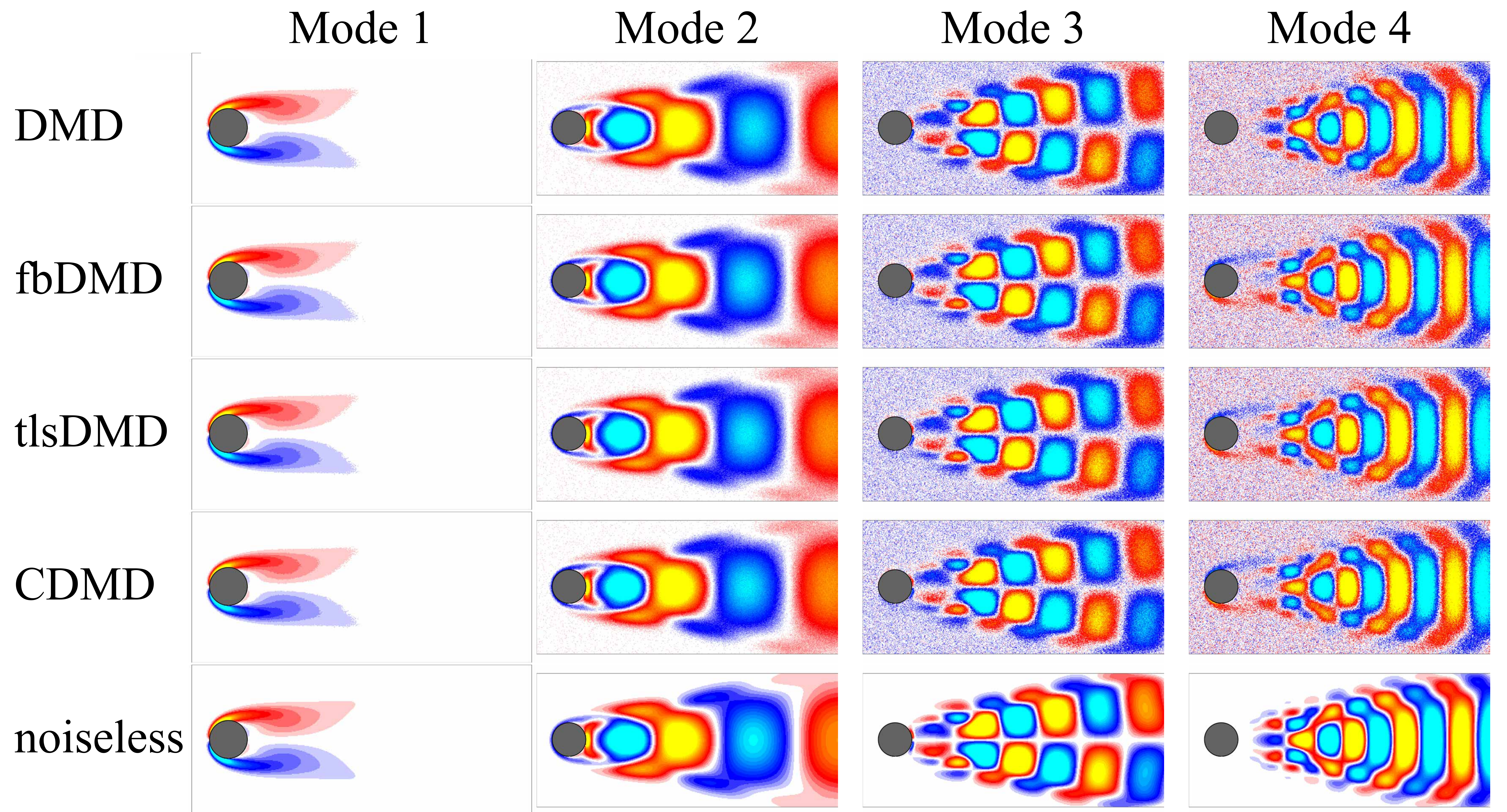}
 \caption{The averaged DMD modes associated with the cylinder wake are shown for data consisting of noise with variance $\sigma^2=0.1$. In spite of the large amounts of white Gaussian noise, all methods produce relatively good estimations when compared to the noiseless scenario (bottom row).}
 \label{fig:cyn_eigvecs}
\end{figure*}

The last example we consider in this work is of a fluid flow past a cylinder simulated using a numerical solver. We obtain a time series of fluid vorticity fields consisting of $n=150$ snapshots regularly sampled in time with $\Delta t = 0.2$. We refer to~\cite{kutz2016dynamic} for additional details regarding this dataset such as the chosen physical parameters and other numerical considerations. It is important to note that this particular flow is inherently \emph{non-linear} and thus the underlying assumptions of methods such as \code{optimized DMD} may not hold. Specifically, it is unclear which functions to fit and whether exponential functions are a good choice in this scenario. In contrast, our approach (as well as other DMD techniques) does not impose restricting conditions on the input data, making it applicable in such challenging scenarios. In~\Cref{fig:cyn_eigvals}, we repeatedly compute the eigenvalues associated with a noisy version of the input data for various noise levels, and we plot the average results as compared to the estimates obtained from the clean observations. Specifically, we repeat this experiment $N=1000$ times for noise with variance $\sigma^2 = 0.001,0.01,0.1$ and $\text{SNR}=30, 20, 10\ \text{dB}$, respectively. Clearly, \code{Exact DMD} exhibits a bias in its estimations which is consistent with previous reports such as~\cite{dawson2016characterizing}. On the other hand, \code{fbDMD} and \code{tlsDMD} generate improved approximations of the eigenvalues with less accuracy as the noise increases. Our approach is successful in measuring nearly zero growth for all eigenvalues and noise levels with a bias in frequencies for the least dominant eigenvalues. In~\Cref{fig:cyn_eigvecs}, we demonstrate the averaged dominant DMD modes obtained for $\sigma^2=0.1$. In this case, all methods perform comparably well in the noiseless case, where the averaged modes associated with less dominant eigenvalues are clearly noisier.

\section{Discussion and Future Work}
\label{sec:conclusions}

In this work, we presented a new method for computing Dynamic Mode Decomposition operators that is based on a variational formulation of the underlying problem, while taking into account the forward and backward dynamics. The obtained minimization is solved using an effective splitting ADMM scheme, which performs well in practice in terms of computational requirements and achieved accuracy. Moreover, it is shown that CDMD could be modified to a provably convergent ADMM scheme at the cost of insignificant additional computations. We demonstrate the performance of our method on a few benchmark dynamical systems, compared to several state-of-the-art approaches. Our conclusion is that the generality of our model, along with its improved accuracy for high levels of noise and low number of observations, makes it an interesting alternative among current existing techniques.

One limitation of our approach is related to the non-linearity and non-convexity of the problem we aim to solve. In particular, it is not clear at this point whether the obtained minimizers are local or global, which is a general challenge in these type of problems, as was also noted in~\cite{askham2018variable}. Another difficulty associated with our work involves the interplay between the chosen value of the penalty parameter $\rho$ and the obtained solutions. While in general our technique is robust to the initial value of $\rho$ due to scheme~\cref{eq:update_rho}, it still affects our results to some extent, as can be seen in~\Cref{fig:sine_consistency}, where for large values of $n$, our consistency error increases. Finally, our algorithm is more computationally demanding compared to the alternatives. However, this is highly dependent on the particular implementation and choice of parameters such as convergence thresholds and thus can be reduced, depending on the particular application at hand.

We believe that formulating DMD in a variational form is important as other regularizers may be considered along with our consistency constraints such as sparsity promoting penalty terms~\cite{jovanovic2014sparsity}. We leave this consideration for future work. Moreover, we would like to explore the relation of our approach to existing techniques such as \code{tlsDMD}. Another interesting direction is to combine the current work with methods that numerically compute an optimal basis~\cite{wynn2013optimal}. The associated problem is extremely challenging as it is of high dimension, non-linear and typically non-convex. We believe that some of the ideas that we presented in this work could be generalized to this case and we plan on pursuing this direction in the future.

\newpage
\appendix
\section{Convexity of $f(\mca{A},\mca{B})$}
\label{app:f_sconvex}

The function $f(\mca{A},\mca{B})$ is $(m,M)$-strongly convex if each of its terms is strongly convex. Thus, we show it for the first term $\hat{f}(A') = \frac{1}{2}|A'X-Y|_F^2$, and we note that a similar derivation could be carried for the other term. We recall the gradient of $\hat{f}(A')$ and we vectorize it to arrive at the following formulation
\begin{align*}
\nabla \hat{f}(A') = (A'X-Y)X^T = A' X X^T - Y X^T \equiv (X X^T \otimes I) \vecm(A') - \vecm(Y X^T) \ .
\end{align*}
Therefore, when viewed as a vectorized function, the Hessian of $\hat{f}$ is given by $\nabla^2\hat{f} = X X^T \otimes I$. The matrix $X \in \mbb{R}^{r \times n}$ can be assumed to have full rank, since $r \ll n$, and thus $X X^T$ is positive definite (PD). It is known that the product of two PD matrices is also PD, which means that there exists a scalar $m>0$ such that the Hessian $\nabla^2\hat{f}-mI$ is positive semi-definite, and we conclude that $\hat{f}$ is an $m$-strongly convex function. Finally, $\hat{f}$ is also $M$-Lipschitz differentiable since $|(A_1'-A_2')XX^T|_F \leq |XX^T|_F \cdot |(A_1'-A_2')|_F$ and $|XX^T|_F$ is positive and bounded.

\bibliographystyle{siamplain}
\bibliography{cdmd}

\end{document}